\documentclass[letterpaper,10pt]{article}
\usepackage[utf8]{inputenc}
\usepackage{amsmath, amsthm, amsfonts, amssymb}
\usepackage{microtype}
\usepackage{tikz}
\usetikzlibrary{calc}
\title{Non-injectivity of Nonzero Discriminant Polynomials and Applications to Packing Polynomials}
\author{K{\aa}re Schou Gjaldb{\ae}k}

\newtheorem{thm}{Theorem}
\newtheorem{lem}[thm]{Lemma}
\newtheorem{cor}[thm]{Corollary}
\theoremstyle{remark}
\newtheorem*{rem*}{Remark}

\newcommand{\R}{\mathbb{R}} 
\newcommand{\N}{\mathbb{N}} 
\newcommand{\Z}{\mathbb{Z}} 
\newcommand{\cC}{\mathcal{C}} 
\newcommand{\pr}{^{\prime}} 
\newcommand{\il}{\textrm} 
\newcommand{\size}[1]{\#{#1}} 
\newcommand{\area}{\operatorname{area}} 
\newcommand{\bigoh}{\textsl{O}} 
\newcommand{\forwhich}{~:~} 

\newcommand{\picqpp}[6][scale=1.0]{
    \begin{tikzpicture}[#1]
    	\pgfmathsetmacro{\xmin}{-1}
    	\pgfmathsetmacro{\ymin}{-1}
    	\pgfmathtruncatemacro{\n}{#2}
    	\pgfmathtruncatemacro{\m}{#3}
    	\pgfmathtruncatemacro{\k}{#4}
    	\pgfmathsetmacro{\xmax}{#5}
    	\pgfmathsetmacro{\ymax}{#6}
    	\pgfmathtruncatemacro{\l}{gcd(\m-1,\n)}
    	
    	\draw[black!20] ({0.1 + \xmin}, {0.1 + \xmin}) grid ({\xmax - 0.1}, {\ymax - 0.1});
    	\begin{scope}
    		\clip (\xmin,\ymin) rectangle (\xmax + 0.1,\ymax);
    		\pgfmathtruncatemacro{\s}{\n/\l};
    		\pgfmathtruncatemacro{\stairs}{\xmax*\s + (\m < 1)*(\xmax*\s + 3)};
			\pgfmathtruncatemacro{\firstchunk}{\ymax/\l + 1}
			\pgfmathtruncatemacro{\secondchunk}{\xmax * \n/\l - \ymax*(\m-1)/\l + 1}
			\pgfmathtruncatemacro{\thirdchunk}{\xmax * \s + 1}
    		\foreach \i in {0, ..., \stairs}{
				\ifnum \i<\firstchunk {
					\draw[black!60] ({\i/\s},0)--({\m*\i/\s}, {\n*\i/\s});
				}\else{
					\ifnum \i<\secondchunk {
						\draw[black!60] ({\i/\s},0)--({(\m-1)*\ymax/\n + \i/\s}, \ymax);
					}\else{
						\ifnum \i<\thirdchunk {
							\draw[black!60] ({\i/\s},0)--(\xmax, {\n*\xmax/(\m - 1) - \l*\i/(\m - 1)});
						}\else{
							\ifnum \m=0{
								\draw[black!60] ({\i/\s},0)--(\xmax, {\n*\xmax/(\m - 1) - \l*\i/(\m - 1)});
							}\else{
							}\fi
						}\fi
					}\fi
				} \fi
			}
			\ifnum \m>0{
				\draw[very thick, black] (0,0) -- (\xmax , \xmax*\n/\m);
			}\else{
				\draw[very thick, black] (0,0) -- (0 , \ymax);
			}\fi
     	\end{scope}
    	\draw[thick,-latex] (\xmin,0)--(\xmax + 0.7,0);
    	\draw[thick,-latex] (0,\ymin)--(0,\ymax + 0.3);
    	\foreach \x in {0, ..., \xmax}{
    		\foreach \y in {0, ..., \ymax}{
				\pgfmathtruncatemacro{\pxy}{
					\n/2 * (\x - (\m-1)*\y/\n) * (\x - (\m-1)*\y/\n - \k*\l/\n)
					+ \x
					+ (\k*\l + 1 - \m)*\y/\n
					+ abs(\k) - 1
					+ 0.5
				}
				\ifnum \m=0{
					\node[color=blue!50!black] at (\x,\y) {$\pxy$};
				}\else{
					\pgfmathtruncatemacro{\ycoord}{\x*\n/\m};
					\ifnum \y>\ycoord{
						\node[color=red!80!white] at (\x,\y) {$\pxy$};
					}\else{
						\node[color=blue!50!black] at (\x,\y) {$\pxy$};
					}\fi
				}\fi
    		}
    	}
    \end{tikzpicture}
}

\begin{document}

\maketitle

\begin{abstract}
Define the sector $S(\alpha) := \{ (x,y) \in \R^2 \forwhich 0 \leq y \leq \alpha x \}$.
The sector is called irrational if $\alpha$ is irrational.
A packing polynomial on $S(\alpha)$ is a polynomial which bijectively maps integer lattice points
of $S(\alpha)$ onto the non-negative integers.
We show that an integer-valued quadratic polynomial on $\R^2$
can not be injective on the integer lattice points of any affine convex cone if its discriminant is nonzero.
A consequence is the non-existence of quadratic packing polynomials
on irrational sectors of $\R^2$.
\end{abstract}

\section{Background}
\begin{figure}[hbtp]
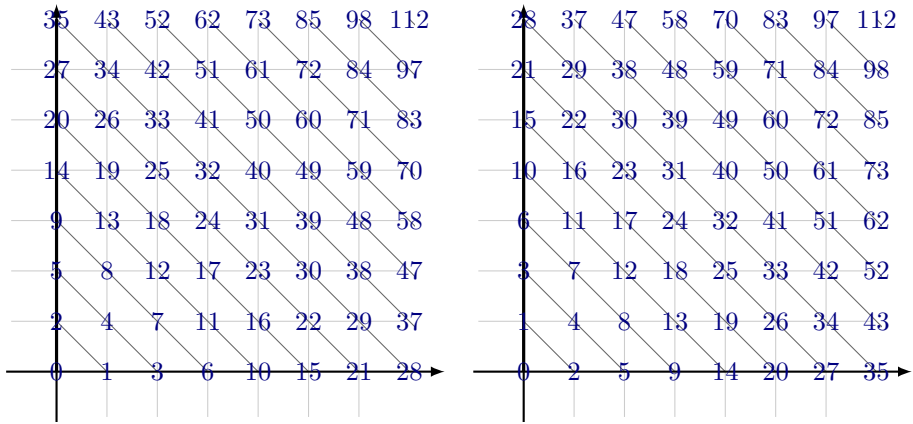

	\label{fig:CantorPols}
	\centering
	\picqpp[scale=0.67]{1}{0}{1}{7}{7}
	\picqpp[scale=0.67]{1}{0}{-1}{7}{7}
	\caption{The two Cantor polynomials correspond to enumerating the lattice along diagonals.}
\end{figure}

	In the seminal 1878 paper \emph{Ein Beitrag zur Mannigfaltigkeitslehre} \cite{cantor1878beitrag},
	Cantor introduces the polynomial
	\[
		f(x,y) = x + \frac{(x + y - 1)(x + y + 2)}{2}
	\]
	which bijectively maps $\N \times \N$ onto $\N$,
	$\N$ denoting the positive integers.
	Translating, and interchanging the variables, leads to the two
	\emph{Cantor Polynomials}
	\begin{align*}
		F(x,y) &= \frac{1}{2}(x+y)(x+y+1) + x\;,\\
		G(x,y) &= \frac{1}{2}(x+y)(x+y+1) + y\;.
	\end{align*}
	These two quadratics are bijections from $\N_0 \times \N_0$ onto $\N_0$,
	$\N_0$ denoting the non-negative integers.
	See Figure \ref{fig:CantorPols}

	In 1923, Fueter and P{\'o}lya \cite{fuetPol23}  prove that the two
	Cantor Polynomials are the only two \emph{quadratic} polynomials
	admitting such a bijection. They further conjecture that a bijection is
	impossible for a polynomial of degree other than two.
	Fueter and P{\'o}lya's proof relies on a corollary of the Lindemann-Weierstra{\ss}
	theorem which says that the trigonometric and hyperbolic functions
	evaluated at non-trivial algebraic numbers are transcendental.
	In 2001, Vsemirnov \cite{vsem01} proves the theorem using elementary methods.
	Theorem \ref{thm:noninj} of this paper allows for Fueter and P{\'o}lya's strategy
	to work without the need for the theorem of Lindemann-Weierstra{\ss}
	(Fueter and P{\'o}lya's result is a special case of
	the classification of rational sectors given in \cite{brandtgjald20}).
	
	In two papers \cite{lewRos78, lewRos78p2} from 1977,
	Lew and Rosenberg develop a more general theory,
	and some of the terminology they introduce has taken hold.
	They provide a partial result on Fueter and P{\'o}lya's conjecture in
	proving the non-existence of polynomials of degree 3 and 4.
	The general problem remains open.

	For some preliminary results,
	Lew and Rosenberg consider polynomials on arbitrary sectors,
	regions that are the convex hull of two half-lines starting at the origin,
	yet they do not study polynomial bijections from general sectors onto $\N_0$.
	This is the subject of Nathanson's 2014 paper
	\emph{Cantor Polynomials for Semigroup Sectors}
	\cite{nath14}.
	Nathanson looks at
	quadratic polynomial bijections from arbitrary sectors in the first quadrant,
	in particular the lattice points in the convex cone spanned by the $x$-axis
	and the line $y = \alpha x$ for some $\alpha > 0$.
	Nathanson determines all such polynomials for sectors given by $\alpha = 1/n$,
	$n \in \N$. Furthermore, he finds two quadratic polynomials for each sector with integer
	$\alpha$-value.
	He raises the question of which rational values for $\alpha$ allow
	for quadratic polynomial bijections, and whether such is possible for
	irrational $\alpha$.
	
	The same year, Stanton and Brandt answer the questions regarding rational sectors.
	In \cite{stan14},
	Stanton determines all quadratic polynomials for $\alpha \in \N$.
	In addition to the polynomials discovered by Nathanson,
	she finds two quadratic packing polynomials for the sectors given by $\alpha = 3,4$
	and proves that there are no more.
	The non-integral rational case is addressed by Brandt
	in the preprint paper \cite{brandt14}.
	
	Regarding the irrational case, to be accurate,
	Nathanson asks the reader to show that no bijections are possible,
	although he never names it a conjecture.
	Corollary \ref{cor:noirr} of this paper obliges.
	
\section{Quadratic Packing Polynomials on Sectors}

	Let $\alpha > 0$. Define the sector
	\[
		S(\alpha) = \{(x,y)\in\R^2 \forwhich 0 \leq y \leq \alpha x\}\;.
	\]
	We let $S(\infty)$ denote the first quadrant, i.e.~the sector pertaining to the original problem.
	Adopting the terminology of Lew and Rosenberg,
	we refer to a bijection from $S(\alpha) \cap \Z^2$ onto $\N_0$ as a \emph{packing function},
	or, in the case of a polynomial, which is our sole focus, \emph{packing polynomial}.

	An immediate prerequisite for a packing polynomial is that it is \emph{integer-valued}.
	That is, it must take integer values on integer lattice points.
	A consequence of standard results on integer-valued polynomials (see e.g.~\cite{lang2002algebra},
	Chp.~X, Lem.~6.4)
	is that a quadratic packing polynomial on any sector must be of the form
	\begin{equation}\label{eq:intvalform}
		P(x,y) = A \frac{x(x-1)}{2} + Bxy + C\frac{y(y-1)}{2} + Dx + Ey + F
	\end{equation}
	with $A,B,C,D,E,F \in \Z$.
	Lew and Rosenberg make the following observation
	(see \cite{lewRos78}, Prop.~3.4).
	\begin{lem}\label{lem:quadpos}
		Let $P(x,y)$ be a packing polynomial%
\footnote{In fact, it is only required for the polynomial to be an injection into $\N_0$,
a class Lew and Rosenberg call \emph{storing functions}.
}
		on the sector $S(\alpha)$.
		If $(m,n) \in S(\alpha)\setminus\{(0,0)\}$ is an integer lattice point, then the homogeneous
		quadratic part of $P(x,y)$:
		\[
			P_2(x,y) = \frac{A}{2} x^2 + B xy + \frac{C}{2} y^2\;,
		\]
		must take only positive values on the ray $\{(xm, xn), x > 0\}$.
	\end{lem}
	
	\noindent An immediate consequence is that we must have $A > 0$.

\section{Non-injectivity when Discriminant is Zero}

Let $\omega_1, \omega_2 \in \R^2$ with $\omega_1 \neq \omega_2$
and define the closed convex cone
\[
	\cC(\omega_1, \omega_2) = \{u \omega_1 + v \omega_2 \forwhich u, v \geq 0 \}
\]
and for $\omega_0 \in \R^2$ the affine convex cone
\[
	\cC_{\omega_0}(\omega_1, \omega_2) = \cC(\omega_1, \omega_2) + \omega_0\;.
\]

\begin{thm}\label{thm:noninj}
	Let $P: \R^2 \to \R$ be an integer-valued
	quadratic polynomial.
	If the discriminant of $P$ is non-zero,
	then $P$ cannot be injective on the integer lattice points of any affine convex cone.
\end{thm}

\begin{proof}
	Let $P(x,y)$ have the form \eqref{eq:intvalform}.
	We will denote its discriminant $\Delta = B^2-AC$
	and use the shorthand notation $D\pr = D-\frac{A}{2}$, $E\pr = E - \frac{C}{2}$.
	Let $\cC = \cC_{\omega_0}(\omega_1, \omega_2)$ be an arbitrary affine cone.
	Fix coprime integers $r, s$ with $r\neq 0$. Every lattice point lies on a line
	$L^{(i)}_{r,s}: y = \frac{s}{r} x + \frac{i}{r}$ for a unique $i$.
	For each $i$, consider the restriction of $P$ to $L^{(i)}_{r,s}$:
	\begin{align*}
		Q_i(x) &= P\left(x, \frac{s}{r} x + \frac{i}{r}\right)\\
			&= \frac{1}{2 r^2}(A r^2 + 2B r s + C s^2) x^2\\
				&\hspace{1cm}+ \frac{1}{r^2}((B r + C s)i + r(D \pr r + E\pr s)) x
				+ \il{const.}
	\end{align*}
	If $r, s$ are chosen such that $A r^2 + 2B r s + C s^2 \neq 0$,
	then the values of $Q_i(x)$ are symmetric around
	\[
		x_i = -\frac{(B r + C s)i + r (D \pr r + E\pr s)}{A r^2 + 2B r s + C s^2}\;.
	\]
	The corresponding $y$-coordinate on $L^{(i)}_{r,s}$ is
	\[
		y_i = \frac{(Ar + Bs)i - s (D\pr r + E\pr s)}{A r^2 + 2B r s + C s^2}\;.
	\]
	This means that, for any choice of $r,s$ with $A r^2 + 2B r s + C s^2 \neq 0$,
	we have
	\[
		P(x_i + r, y_i + s) = P(x_i - r, y_i -s)
	\]
	for all $i$.
	
	These points of symmetry, $(x_i, y_i)$, fall on the straight line, $L^{\il{sym}}_{r,s}$, with slope
	$-\frac{Ar + Bs}{Br + Cs}$ which passes through the point
	\[
		(x_0, y_0) = \left(\frac{CD\pr - BE\pr}{\Delta}, \frac{AE\pr - BD\pr}{\Delta} \right)\;.
	\]
	(To see this, replace $i$ with $\frac{1}{\Delta} ((AE\pr - BD\pr)r - (CD\pr - BE\pr)s)$
	in the formulas for $x_i, y_i$.)
	Note that the point $(x_0, y_0)$ does not depend on the choice of $r, s$.%
\footnote{
	When $\Delta \neq 0$, $P(x,y)$ can be rewritten as
	\[
		P(x,y) = \frac{A}{2}(x-x_0)^2 + B(x-x_0)(y-y_0) + \frac{C}{2} (y-y_0)^2
			+ \frac{D\pr}{2} x_0 + \frac{E\pr}{2} y_0 + F\;.
	\]
	The point $(x_0, y_0)$ is the center of the level curves $P(x,y) = n$ which are
	ellipses when $\Delta <0$, or hyperbolas when $\Delta > 0$.
}
	We want to choose $r, s$ such that $(x_i + r, y_i + s)$ and $(x_i - r, y_i -s)$
	are both lattice points in $\cC$ for some $i$.
	This will violate injectivity.
	
	Put $\cC_0 = \cC_{(x_0, y_0)}(\omega_1, \omega_2)$ and pick an arbitrary lattice point
	$(m,n) \in \cC \cap \cC_0$. Choosing
	\[
		\frac{s}{r} = - \frac{A(m-x_0) + B(n-y_0)}{B(m-x_0) + C(n-y_0)}
	\]
	in lowest terms, we have $-\frac{Ar + Bs}{Br + Cs} = \frac{n-y_0}{m-x_0}$.
	This means that
	$L^{\il{sym}}_{r,s}$ passes through the lattice point $(m,n) \in \cC \cap \cC_0$
	and therefore infinitely many since its slope is rational.
	So for infinitely many $i$, $(x_i, y_i)$ is a lattice point in
	$\cC \cap \cC_0$ and eventually we will find an $i$ for which
	both $(x_i + r, y_i + s)$ and $(x_i - r, y_i -s)$ are lattice points in $\cC$.
	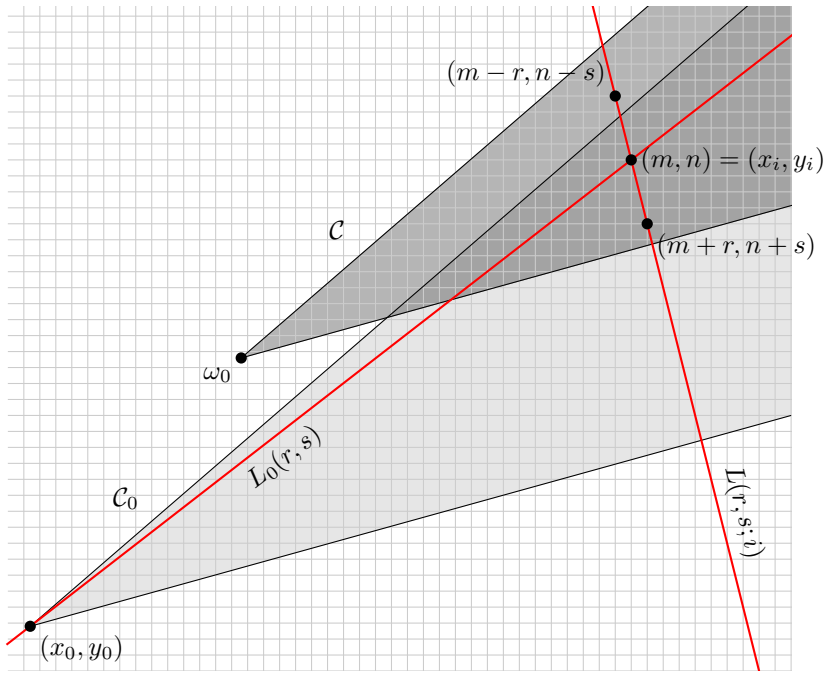
\begin{figure}[htbp]
		\label{fig:cones}
 		\centering
		\begin{tikzpicture}[scale=0.85]
			\def\xmax{12}
			\def\ymax{10}
			\def\xnul{0.1}	\def\ynul{0.7}
			\def\xb{12}		\def\yb{11}
			\def\xc{12}		\def\yc{4}
			\def\xv{3.3}	\def\yv{4.2}
			\newcommand{\coneintcol}{black!36}
			\newcommand{\conecol}{black!29}
			\newcommand{\conezerocol}{black!10}
			\newcommand{\gridcol}{black!20}
			\coordinate (O) at (\xnul, \ynul);
			\coordinate (V) at (\xv, \yv);
			\def\ppp{0.25}
			\def\qqq{-1}
			\coordinate (L) at (9.5,8);
			\coordinate (Lm) at ($(L) - (\ppp,\qqq)$);
			\coordinate (Lp) at ($(L) + (\ppp,\qqq)$);
			\coordinate (A) at (O);
			\coordinate (B) at (\xb, \yb);
			\coordinate (C) at (\xc, \yc);
			\coordinate (D) at ($(A) + (V)$);
			\coordinate (E) at ($(B) + (V)$);
			\coordinate (F) at ($(C) + (V)$);
			\pgfmathsetmacro{\xint}{
				(
					\xnul * (\xc - \xnul) * (\yb - \ynul) -
					(\xnul + \xv) * (\yc - \ynul) * (\xb - \xnul) +
					(\xc - \xnul) * (\xb - \xnul)*(\yv)
				) /
				(
					(\xc - \xnul) * (\yb - \ynul) -
					(\yc - \ynul) * (\xb - \xnul)
				)
			}
			\pgfmathsetmacro{\yint}{
				(\yb - \ynul)*(\xint - \xnul) / (\xb - \xnul)  + \ynul
			}
			\coordinate (I) at (\xint, \yint);
			\def\xstart{-0.4}
			\def\xstop{5.1}
			\begin{scope}
				\clip (-1,0) rectangle (12, 10.4);
				\filldraw[\conezerocol] (A)--(B)--(C)--cycle;
				\filldraw[\conecol] (D)--(E)--(F)--cycle;
				\filldraw[\coneintcol] (I)--(B)--(F)--cycle;
				\foreach \y in {0, ..., \ymax}{
					\draw[\gridcol] (-0.25,\y)--(\xmax + 1 , \y);
					\draw[\gridcol] (-0.25,\y + 0.25)--(\xmax + 1, \y + 0.25);
					\draw[\gridcol] (-0.25,\y + 0.50)--(\xmax + 1, \y + 0.50);
					\draw[\gridcol] (-0.25,\y + 0.75)--(\xmax + 1, \y + 0.75);
				}
				\foreach \x in {0,..., \xmax}{
					\draw[\gridcol] (\x, -0.25)--(\x, \ymax + 1);
					\draw[\gridcol] (\x + 0.25, -0.25)--(\x + 0.25, \ymax +1);
					\draw[\gridcol] (\x + 0.50, -0.25)--(\x + 0.50, \ymax +1);
					\draw[\gridcol] (\x + 0.75, -0.25)--(\x + 0.75, \ymax +1);
				}
				\draw (A)--(B);
				\draw (A)--(C);
				\draw (D)--(E);
				\draw (D)--(F);
				\draw[thick, red, shorten >= -2.9 cm, shorten <= -0.4 cm] (O) -- node[black, below, pos=0.4, sloped] {$L_0(r,s)$} (L);
				\draw[thick, red, shorten >= -6.7 cm, shorten <= -1.8 cm] (Lm)-- node[black, above, pos=3.3, sloped] {$L(r,s; i)$} (Lp);
			\end{scope}
			\node[circle, fill, inner sep=1.5pt] at (O) {};
			\node[anchor = north west] at (O) {$(x_0,y_0)$};
			\node[circle, fill, inner sep=1.5pt] at (D) {};
			\node[anchor = north east] at (D) {$\omega_0$};
			\node[circle, fill, inner sep=1.5pt] at (L) {};
			\node[anchor = west] at (L) {$(m,n) = (x_i, y_i)$};
			\node[circle, fill, inner sep=1.5pt] at (Lm) {};
			\node[anchor = south east] at (Lm) {$(m - r, n - s)$};
			\node[circle, fill, inner sep=1.5pt] at (Lp) {};
			\node[anchor = north west] at (Lp) {$(m + r, n + s)$};
			
			\node at ($ (D) + (1.5,2) $) {$\cC$};
			\node at ($ (O) + (1.5,2) $) {$\cC_0$};
		\end{tikzpicture}

		\caption{
			$L_0(r,s)$ with $\frac{r}{s} = - \frac{A(m-x_0) + B(n-y_0)}{B(m-x_0) + C(n-y_0)}$
			passing through a lattice point in $\cC \cap \cC_0$.
		}
	\end{figure}
\end{proof}

\section{Applications to Packing Polynomials}

Let $P(x,y)$ be a quadratic packing polynomial on the sector $S(\alpha)$ with discriminant $\Delta$.
$S(\alpha)$ is a (affine) convex cone, so by the previous section we must have $\Delta = 0$.
Employing the strategy of Lew and Rosenberg \cite{lewRos78}, we consider the regions
\[
	R_n = \{(x,y) \in S(\alpha) \forwhich 0 \leq P(x,y) \leq n \}\;.
\]
For a packing polynomial, each region $R_n$ contains $n + 1$ lattice points.
This means that it is a necessary condition that
\begin{equation}\label{eq:limitisone}
	\lim_{n \to \infty} \left( \frac{1}{n} \size{(R_n \cap \Z^2)}\right) = 1 \;.
\end{equation}
Furthermore, we have
\begin{lem}\label{lem:errorterm}
	If $P(x,y)$ is a quadratic packing polynomial on $S(\alpha)$, then
	\[
		\size{(R_n \cap \Z^2)} = \area(R_n) + \bigoh(\sqrt{n})\;,
	\]
	as $n \to \infty$.
\end{lem}
\begin{proof}
	Since the discriminant is zero,
	the level curves $P(x,y) = n$ are parabolas
	and either $R_n$ is bounded for all $n$ or all level curves,
	including for negative $n$,
	fall inside $S(\alpha)$ which is impossible if $P$ is a packing polynomial.
	We can therefore apply a theorem of Davenport (\cite{dav51}, p.~180) to estimate the number
	of lattice points in each region $R_n$. We have
	\[
		|\area(R_n) - \size{(R_n \cap \Z^2)}| \leq h( |\pi_x (R_n)| + |\pi_y (R_n)|) + h^2 \;,
	\]
	where $|\pi_x (R_n)|$ and $|\pi_y (R_n)|$ denotes the lengths of the projections onto the
	$x$- and $y$-axis, respectively, and $h$ is a fixed constant.%
\footnote{
	The constant $h$ denotes the maximum number of disjoint intervals one can obtain from intersecting
	$R_n$ with a line parallel to one of the coordinate axes.
	Our only concern is that it is bounded by some value.
}
	The length of $y$-projection is bounded by the level curves intersection with the line
	$y = \alpha x$ or the topmost point on the parabola.
	The length of $x$-projection is bounded by the level curves intersection with the
	$x$-axis, the line $y = \alpha x$ or the rightmost point on the parabola.
	Either is $\bigoh(\sqrt{n})$ as $n \to \infty$.
\end{proof}

At this point, we want to note that, since $B^2 = AC$,
the homogeneous quadratic part of $P(x,y)$,
\[
	P_2(x,y) = \frac{1}{2}(\sqrt{A} x \pm \sqrt{C} y)^2 \;,
\]
is non-negative and vanishes only on the rational line
$Ax + By = 0$.
By Lem.~\ref{lem:quadpos}, this means that $P_2(x,y)$ is strictly
positive inside $S(\alpha)$ except at the origin.

\begin{lem}\label{lem:limit}
	If $P(x,y)$ is a quadratic packing polynomial on $S(\alpha)$, then
	\[
		\area(R_n) = \frac{n}{2} \int_0^{\arctan\alpha}
			\frac{d\theta}{\frac{A}{2}\cos^2\theta + B \cos\theta\sin\theta + \frac{C}{2}\sin^2\theta}
			+\bigoh(\sqrt{n})
	\]
	as $n \to \infty$.
\end{lem}

\begin{proof}
	Switching to polar coordinates, the equations of the level curves take the form
	\begin{align*}
		r^2 p_2(\theta) + r p_1(\theta) + F = n \;,
	\end{align*}
	where
	$p_2(\theta) = \frac{A}{2}\cos^2\theta +  B \cos\theta\sin\theta + \frac{C}{2} \sin^2\theta$ and
	$p_1(\theta) = D\pr \cos\theta + E\pr \sin\theta$.
	So $r = \bigoh(\sqrt{n})$ and $r^2 = \frac{n}{p_2(\theta)} + \bigoh(\sqrt{n})$.
	If $A_0$ denotes the (possibly empty) area of $S(\alpha)$ bounded by the level curve $P(x,y) = 0$,
	then the area of $R_n$ is given by
	\begin{align*}
		\area(R_n) &= \frac{1}{2} \int_{0}^{\arctan \alpha} r^2d \theta - A_0
			= \frac{1}{2} \int_{0}^{\arctan \alpha} \frac{n}{p_2(\theta)}  d \theta + \bigoh(\sqrt{n})\;.
	\end{align*}
\end{proof}

\begin{thm}\label{thm:alpharestriction}
	If $P(x,y)$ is a quadratic packing polynomial on $S(\alpha)$, then
	\[
		\alpha = \frac{A}{1-B}\;.
	\]
\end{thm}
\begin{proof}
	Using Lem.~\ref{lem:errorterm} and \ref{lem:limit}, we can calculate the limit
	from \eqref{eq:limitisone} by computing the integral
	\begin{align*}
		\lim_{n \to \infty} \left( \frac{1}{n} \size{(R_n \cap \Z^2)}\right)
			&= \frac{1}{2} \int_0^{\arctan\alpha}
				\frac{d\theta}{\frac{A}{2}\cos^2\theta + B \cos\theta\sin\theta + \frac{C}{2}\sin^2\theta}\\
			&= \int_0^{\alpha}\frac{dt}{A + 2Bt + Ct^2}\;,
	\end{align*}
	applying the variable change $t = \tan\theta$.
	Since, by Thm.~\ref{thm:noninj}, $B^2 = AC$, we either have $B=C=0$, in which case
	\[
		\int_0^{\alpha}\frac{dt}{A + 2Bt + Ct^2} = \int_0^{\alpha}\frac{dt}{A} = \frac{\alpha}{A}\;,
	\]
	or
	\begin{align*}
		\int_0^{\alpha}\frac{dt}{A + 2Bt + Ct^2}
			&= \frac{1}{C} \int_0^{\alpha}\frac{d\theta}{\left(t + \frac{B}{C}\right)^2}
			= \frac{1}{B} - \frac{1}{\alpha C + B}\;.
	\end{align*}
	As noted above, this must be equal to $1$ if $P$ is a packing polynomial.
	Solving for $\alpha$, we get the desired result.
\end{proof}

An immediate consequence of Thm.~\ref{thm:alpharestriction} is the following.
\begin{cor}\label{cor:noirr}
	There are no irrational sectors allowing for quadratic packing polynomials.
\end{cor}

\bibliographystyle{abbrv}
\bibliography{refs}

\begin{thebibliography}{10}

\bibitem{brandt14}
M.~Brandt.
\newblock Quadratic packing polynomials on sectors of {$\R^2$}.
\newblock {\em arXiv:1409.0063v1}, 2014.

\bibitem{brandtgjald20}
M.~Brandt and K.~Gjaldb{\ae}k.
\newblock Classification of quadratic packing polynomials on sectors of
  {$\R^2$}.
\newblock {\em In preparation}, 2021.

\bibitem{cantor1878beitrag}
G.~Cantor.
\newblock Ein {B}eitrag zur {M}annigfaltigkeitslehre.
\newblock {\em Journal fur die reine und angewandte Mathematik}, 84:242--258,
  1878.

\bibitem{dav51}
H.~Davenport.
\newblock On a principle of {L}ipschitz.
\newblock {\em Journal of the London Mathematical Society}, 1(3):179--183,
  1951.

\bibitem{fuetPol23}
R.~Fueter and G.~P{\'o}lya.
\newblock Rationale {A}bz{\"a}hlung der {G}itterpunkte.
\newblock {\em Vierteljschr. Naturforsch. Ges. Z{\"u}rich}, 58:380--386, 1923.

\bibitem{lang2002algebra}
S.~Lang.
\newblock {\em Algebra}.
\newblock Springer, 2002.

\bibitem{lewRos78}
J.~S. Lew and A.~L. Rosenberg.
\newblock Polynomial indexing of integer lattice-points {I}. {G}eneral concepts
  and quadratic polynomials.
\newblock {\em Journal of Number Theory}, 10(2):192--214, 1978.

\bibitem{lewRos78p2}
J.~S. Lew and A.~L. Rosenberg.
\newblock Polynomial indexing of integer lattice-points {II}. {N}onexistence
  results for higher-degree polynomials.
\newblock {\em Journal of Number Theory}, 10(2):215--243, 1978.

\bibitem{nath14}
M.~B. Nathanson.
\newblock Cantor polynomials for semigroup sectors.
\newblock {\em Journal of Algebra and its Applications}, 13(5), 2014.

\bibitem{stan14}
C.~Stanton.
\newblock Packing polynomials on sectors of {$\R^2$}.
\newblock {\em Integers}, 14, 2014.

\bibitem{vsem01}
M.~A. Vsemirnov.
\newblock Two elementary proofs of the {F}ueter-{P}{\'o}lya theorem on pairing
  polynomials.
\newblock {\em Algebra i Analiz}, 13(5):1--15, 2001.

\end{thebibliography}

\end{document}